\documentclass[a4paper]{amsart}

\usepackage[pdftex]{graphicx}
\usepackage{amsmath, amsfonts, amssymb, amsthm}

\usepackage{overcite}
\makeatletter
\def\@cite#1{$\m@th^{\hbox{\@ove@rcfont#1)}}$}
\renewcommand\@biblabel[1]{#1)}
\makeatother

\title[Numerical analysis of transient orbits]{Numerical analysis of transient orbits by the pullback method for covariant Lyapunov vector}
\author{Takayuki Yamaguchi}
\address{Department of Mathematics, Graduate School of Science, Hiroshima University, 1-3-1 Kagamiyama, Higashi-Hiroshima, Hiroshima, 739-8526, Japan}
\email{ytk@hiroshima-u.ac.jp}

\author{Makoto Iima}
\address{Department of Mathematical and Life Sciences, Graduate School of Science, Hiroshima University, 1-3-1 Kagamiyama, Higashi-Hiroshima, Hiroshima, 739-8526, Japan}
\email{makoto@mis.hiroshima-u.ac.jp}

\begin{document}
\begin{abstract}
In order to analyze structure of tangent spaces of a transient orbit,
we propose a new algorithm which pulls back vectors in tangent spaces along the orbit
by using a calculation method of covariant Lyapunov vectors.
As an example, the calculation algorithm has been applied to a transient orbit converging to an equilibrium
in a three-dimensional ordinary differential equations.
We obtain vectors in tangent spaces that converge to eigenvectors of the linearized system at the equilibrium.
Further, we demonstrate that an appropriate perturbation calculated by the vectors can lead
an orbit going in the direction of an eigenvector of the linearized system at the equilibrium.
\end{abstract}

\maketitle

\section{Introduction}

Covariant Lyapunov vectors are vectors on the tangent space at each point on an orbit
and mean directions of perturbations to the orbit
whose growth rates are described by the corresponding Lyapunov exponents.
Practical algorithms to calculate numerically covariant Lyapunov vectors have been developed
\cite{GinelliPoggiEtAl2007,WolfeSamelson2007,KuptsovParlitz2012}
and have been applied to analyze various systems.
Unlike such Lyapunov vectors that are referred to as Gram-Schmidt vectors,
the covariant Lyapunov vectors have an important feature
that their evolution in both forward and backward time also give covariant Lyapunov vectors.

Covariant Lyapunov vectors are defined as vectors of Oseledec splitting,
which is a set of intersections of two Oseledec subspaces.
The method of Ginelli et al., on which we focused on in this paper,
does not calculate explicitly the intersection of two subspaces.
The method calculates repeatedly pullbacks of vectors in tangent spaces along an orbit
after it calculates Gram-Schmidt vectors
and then obtains convergences to covariant Lyapunov vectors of the orbit.

Because there exist no Oseledec splitting for transient orbits
we can not define covariant Lyapunov vectors of transient orbits.
In this paper, however, we are interested in tangent spaces of points
on a transient orbit converging to an equilibrium.
We propose an application of the method of Ginelli et al.
to transient orbits even though covariant Lyapunov vectors are not defined.
Of course, our obtained vectors for transient orbits are not covariant Lyapunov vectors,
but they can be used to characterize the tangent spaces of points on the transient orbit,
because they are vectors that converge to covariant Lyapunov vectors of the equilibrium,
that is, eigenvectors of the linearized system at the equilibrium.

Based on this idea, 
we demonstrate an application of the method of Ginelli et al.
in a simple three-dimensional ordinary differential equation.
In particular, a perturbation of the transient orbit in the direction of the obtained vector,
can create a new orbit that passes quite near the direction of the eigenvector
of the linearized system in a neighborhood of the equilibrium.

In the following sections, we show outline of the calculation algorithm and
a simple example of the application.
We end our paper with a conclusion.

\section{Calculation algorithm}

Our calculation algorithm is a direct application of the method of Ginelli et al \cite{GinelliPoggiEtAl2007}.
We consider an $n$-dimensional ordinary differential equation
\begin{align}
  \dot{u} = g(u, t).
\end{align}
The system of infinitesimal perturbations is described by the following ordinary differential equation;
\begin{align}
  \dot{v} = J(u, t) v, \label{ode-tangent-space}
\end{align}
where $J(u, t)$ is the Jacobian matrix $\frac{\partial g}{\partial u}$.

We assume that an orbit converging to an equilibrium has been obtained.
We select suitable algorithm for time evolution of the orbit.

Along the orbit, we repeat evolving $n$-dimensional volumes in tangent spaces and
normalizing them by QR decomposition.
For a matrix whose columns are vectors on a tangent space,
we evolve the matrix by (\ref{ode-tangent-space}); that is,
we consider the following ordinary differential equation for a matrix $M(t)$;
\begin{align}
  \dot{M}(t) = J(u, t) M(t). \label{ode-tangent-space-matrix}
\end{align}
Let $\mathcal{F}(t_n, \Delta t) M$
be time evolution of $M$ by (\ref{ode-tangent-space-matrix}) from $t_n$ to $t_{n+1} = t_n + \Delta t$
for small time interval $\Delta t$.
For an initial point $u(0)$ and a suitable initial orthogonal matrix $Q_0$,
we calculate orthogonal matrices $\{ Q_n \}$ and upper triangular matrices $\{ R_n \}$
so that
\begin{align}
  \mathcal{F}(t_n, \Delta t) Q_n = Q_{n+1} R_{n+1},
\end{align}
where $t_n = n \Delta t$.
This procedure continues until $n$ equals to $N_2 \gg N_1$ so that
$u(t_{N_1})$ is sufficiently close to the equilibrium.
Then, the columns of $Q_{n}\; (n>N_1)$ converge reasonably
to Gram-Schmidt vectors at the equilibrium.

We let $C_{N_2}$ be a randomly generated upper triangle matrix and repeat pullback $C_n$ by $R_n$.
Matrices $C_{n - 1}$ for $n = N_2, \dots, 1$ are calculated by
\begin{align}
  \tilde{C}_{n - 1} &= R_n^{-1} C_n, \\
  C_{n - 1} &= \tilde{C}_{n - 1} D_{n - 1}^{-1},
\end{align}
where $D_{n - 1}$ are diagonal matrices whose diagonal elements are norms of column vectors of $\tilde{C}_{n - 1}$.
Because $N_1$ and $N_2$ are large enough,
the columns of $Q_{N_1}$ can be regarded as covariant Lyapunov vectors of the equilibrium.
The columns of $Q_n \; (n = 0, \dots, N_1)$ are
vectors in tangent spaces at points of the orbit,
which converge to covariant Lyapunov vectors of the equilibria in forward time.
However, the calculated vectors that are in tangent spaces of the transient orbit
are not the covariant Lyapunov vectors
because covariant Lyapunov vectors can not be defined in tangent spaces of transient orbits.

\section{Example of application of algorithm}

\begin{figure}
  \begin{minipage}{0.49\hsize}
    \centering
    \includegraphics[width=1.2\hsize]{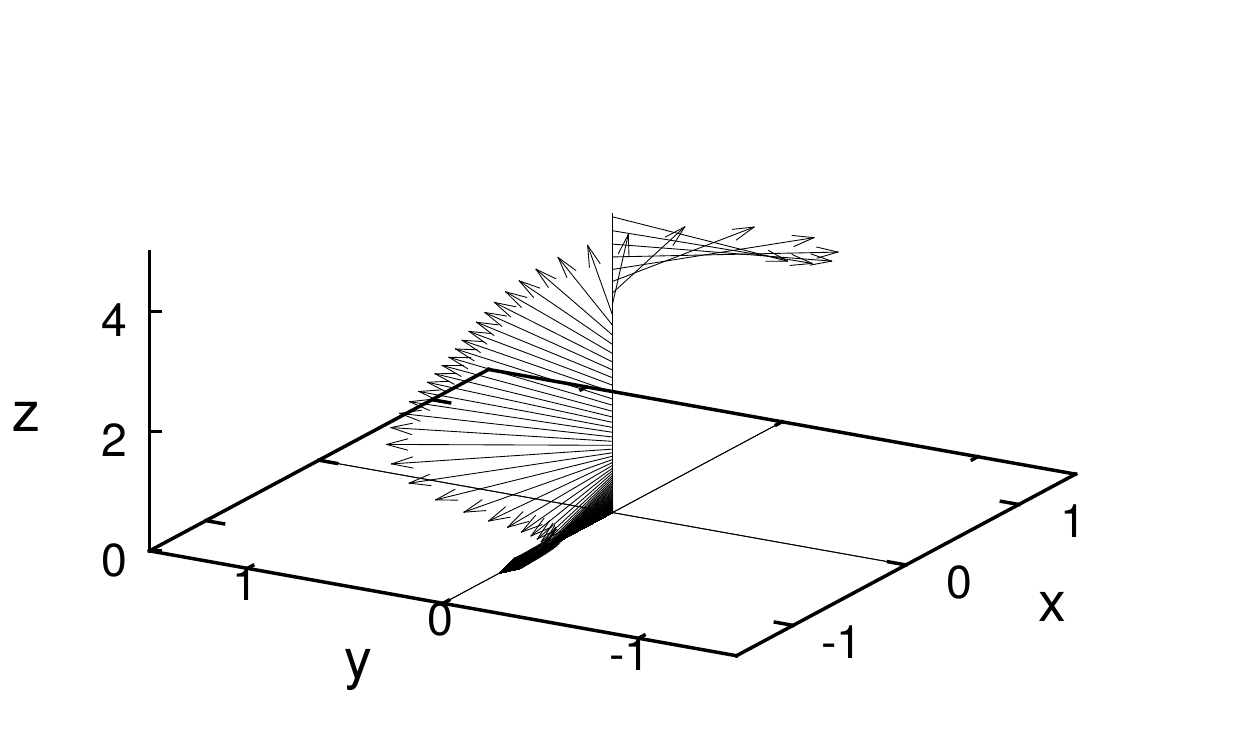}
    Vectors converge to $(-1, 0, 0)^T$ whose eigenvalue is $2$.
  \end{minipage}
  \begin{minipage}{0.49\hsize}
    \centering
    \includegraphics[width=1.2\hsize]{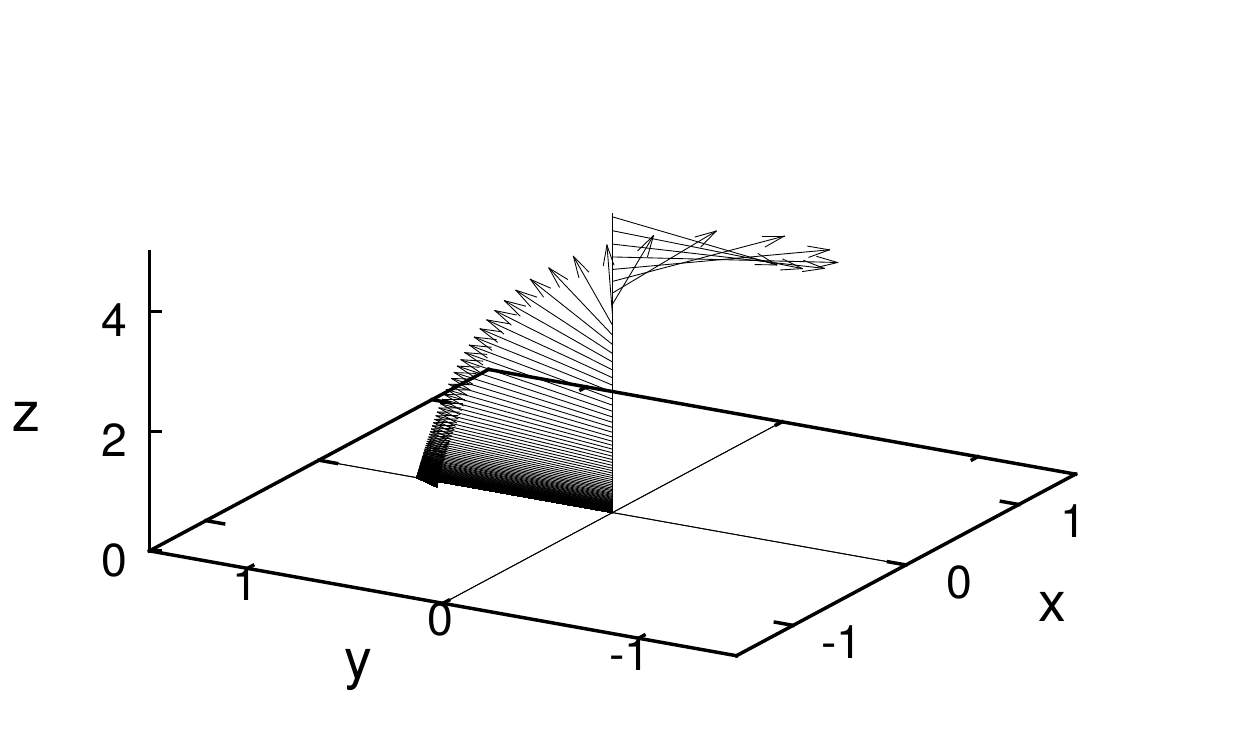}
    Vectors converge to $(0, 1, 0)^T$ whose eigenvalue is $1$.
  \end{minipage}
  \begin{minipage}{0.49\hsize}
    \centering
    \includegraphics[width=1.2\hsize]{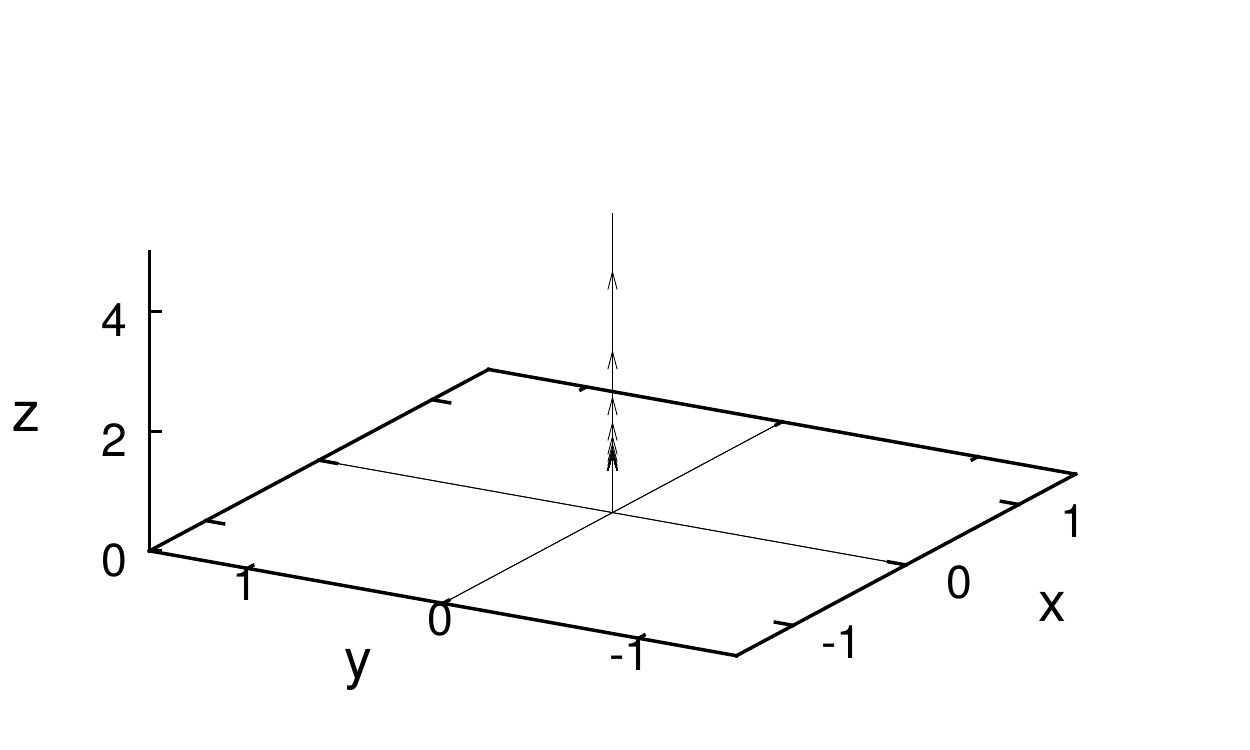}
    Vectors converge to $(0, 0, 1)^T$ whose eigenvalue is $-1$.
  \end{minipage}
  \caption{Vectors converging to eigenvectors of the linearized system at the equilibrium.
  The orbit goes to the equilibrium from $(0, 0, 1000)$.}
  \label{vectors-for-orbit}
\end{figure}

\begin{figure}
  \begin{minipage}{0.49\hsize}
    \includegraphics[width=1.45\hsize]{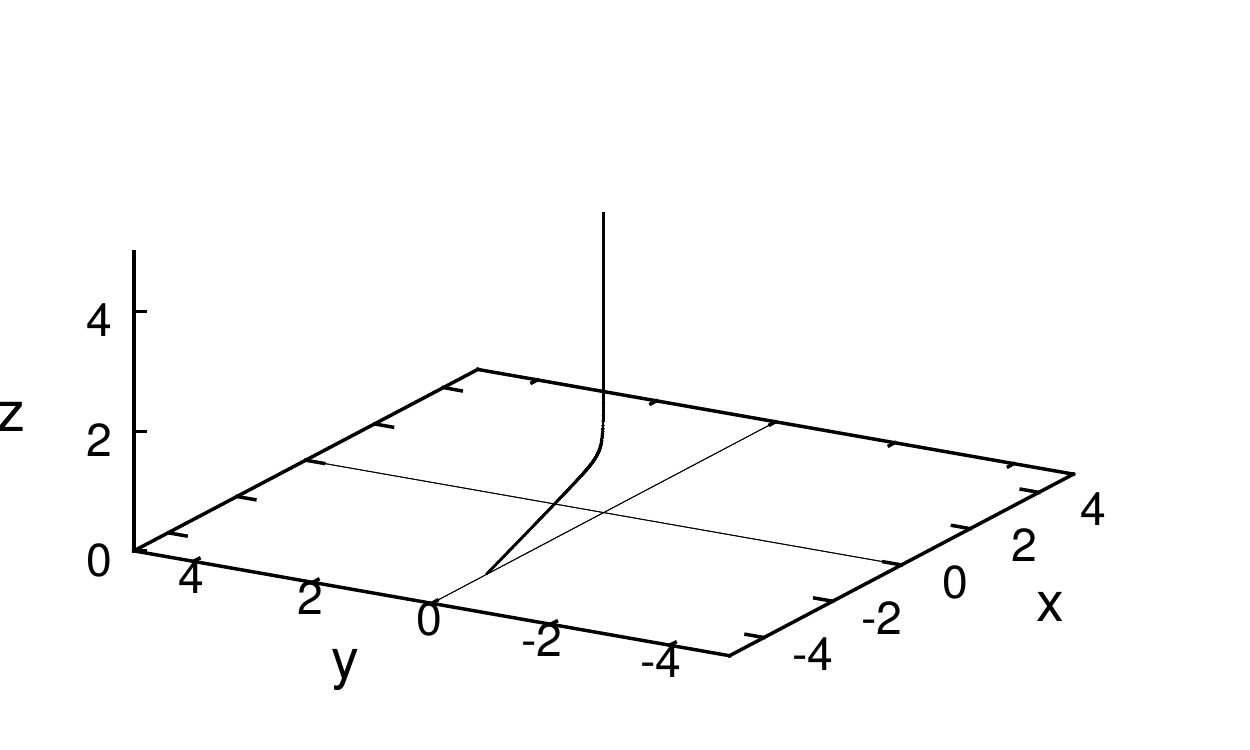}
  \end{minipage}
  \begin{minipage}{0.49\hsize}
    \includegraphics[width=1.25\hsize]{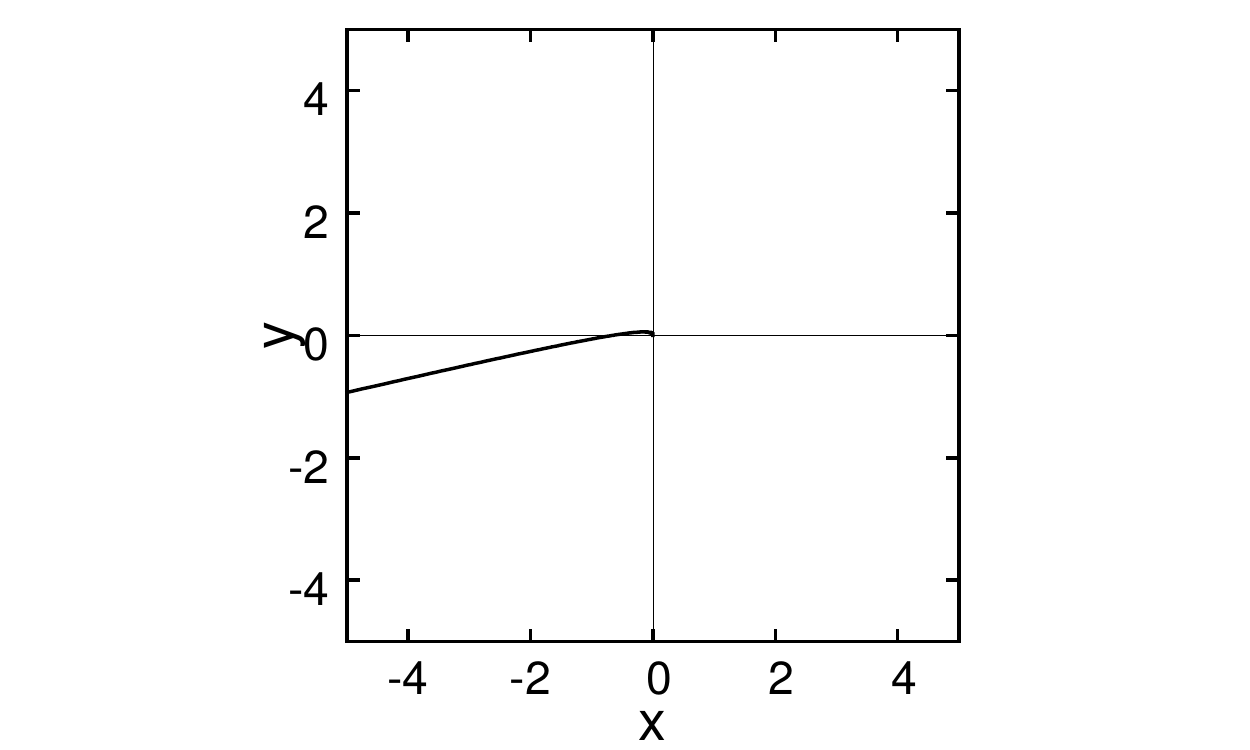}
  \end{minipage}
  \caption{An orbit starting from a point perturbed to the direction of
  the vector converging to the eigenvector whose eigenvalue is $2$.}
  \label{perturbation-of-first-vector}
\end{figure}

\begin{figure}
  \begin{minipage}{0.49\hsize}
    \includegraphics[width=1.45\hsize]{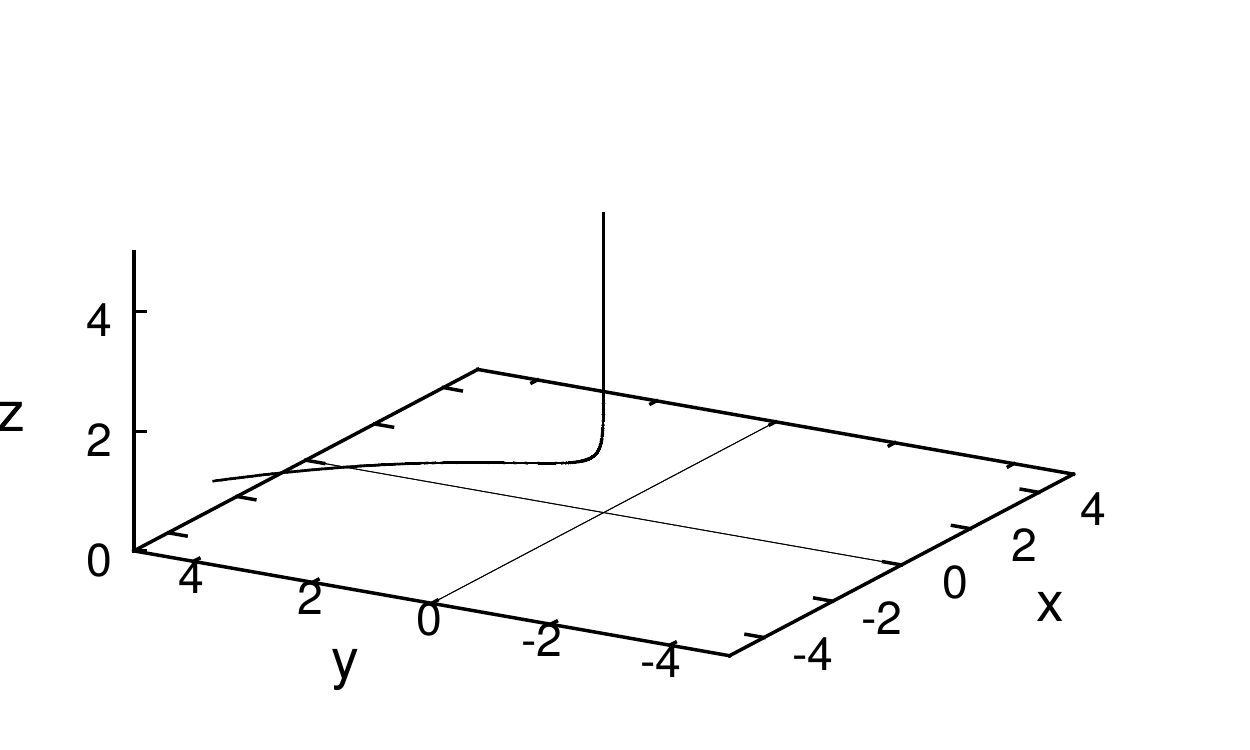}
  \end{minipage}
  \begin{minipage}{0.49\hsize}
    \includegraphics[width=1.25\hsize]{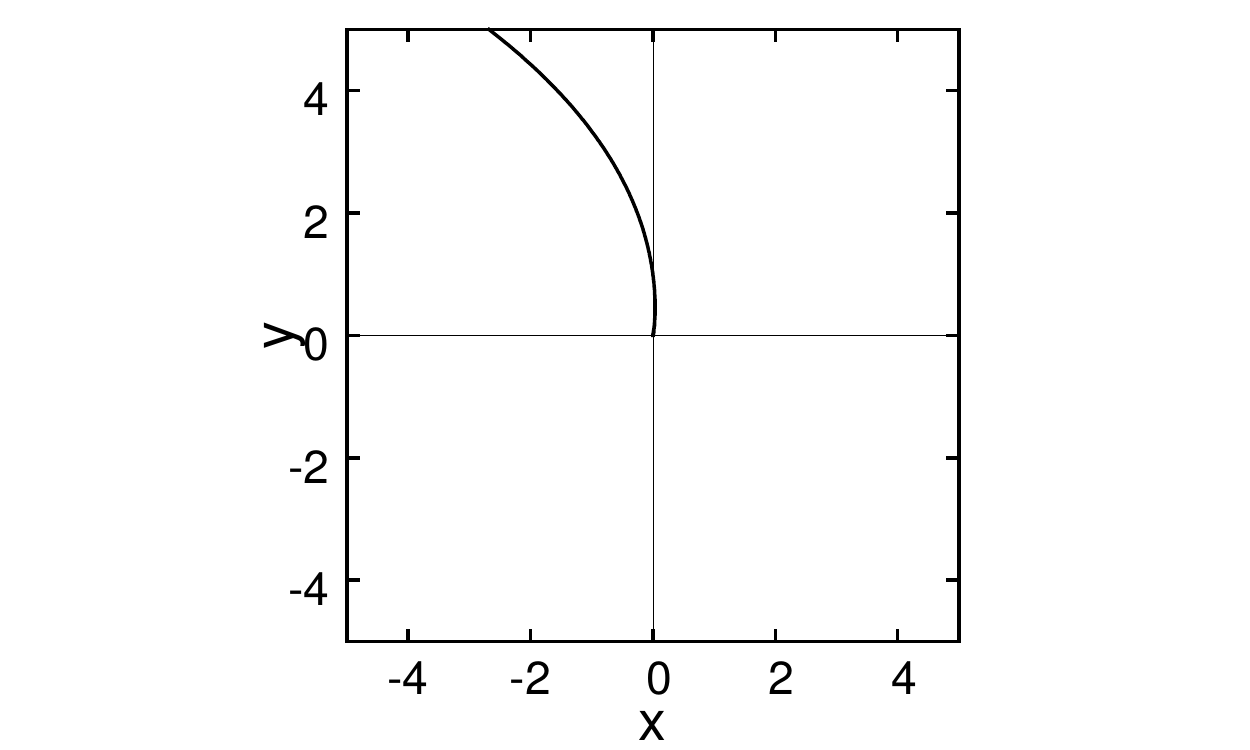}
  \end{minipage}
  \caption{An orbit starting from a point perturbed to the direction of
  the vector converging to the eigenvector whose eigenvalue is $1$.}
  \label{perturbation-of-second-vector}
\end{figure}

To illustrate the ability of the vectors calculated by our method,
we consider three-dimensional ordinary differential equation
$\dot{u} = F(u)$ for $u = (x(t), y(t), z(t))^T$, where
\begin{align}
  F(x, y, z) & = \left(
  \begin{array}{c}
    A x \cos(\arctan(Dz)) - B y \sin(\arctan(Dz)) \\
    A x \sin(\arctan(Dz)) + B y \cos(\arctan(Dz)) \\
    C \arctan(z)
  \end{array}
  \right),
\end{align}
and $A, B, C$ and $D$ are real parameters.
$(x, y, z)^T$ means the transpose of a row vector $(x, y, z)$.
The system has one equilibrium $O = (0, 0, 0)^T$ and
the Jacobian matrix at $O$ is
\begin{align}
  \left(
  \begin{array}{ccc}
    A & 0 & 0 \\
    0 & B & 0 \\
    0 & 0 & C
  \end{array}
  \right).
\end{align}
Clearly, the eigenvalues of the matrix are $A, B, C$ and the corresponding eigenvectors are
$(1, 0, 0)^T$, $(0, 1, 0)^T$, and $(0, 0, 1)^T$.
Orbits starting from $(0, 0, z_0)^T$ for $z_0 > 0$ converges to $O$.

We applied our algorithm to an orbit starting from the point $(0, 0, z_0)$
when the parameter values are $A = 2.0$, $B = 1.0$, $C = -1.0$, and $D = 0.2$.
Then, the eigenvalues of the Jacobian matrix of $F$ at $O$ are $2$, $1$, and $-1$.
The linearized system at $O$ has two unstable directions $(1, 0, 0)^T$ and $(0, 1, 0)^T$ and
one stable direction $(0, 0, 1)^T$.
We have chosen $z_0=1000$ and $\Delta t = 0.1$, $t_{N_1} = 1500$, and $t_{N_2} = 3000$.

Fig. \ref{vectors-for-orbit} shows the vectors obtained by our calculation.
Each of the vectors in the three figures converges to $(-1, 0, 0)^T$, $(0, 1, 0)^T$, and $(0, 0, 1)^T$,
the eigenvectors of the linearized system at $O$, respectively.

We remark that the vectors in Fig. \ref{vectors-for-orbit} have
information about the time evolution of perturbations added to the transient orbit.
Let us consider the orbits starting from initial points in a neighborhood of $(0, 0, z)$.
The orbits go near z-axis and then bend sharply in a direction parallel to xy plane
in a neighborhood of $O$,
because z-axis is the stable manifold of $O$ and
xy plane is the unstable manifold of $O$.
It is difficult to find such an orbit going in the second unstable direction: $(0, 1, 0)^T$,
because we have two unstable direction of $O$
and most orbits go in the largest unstable direction.
However, using the vectors constructed by our method,
we can select such orbit that passes near the second unstable direction.

Let $v_1$ and $v_2$ be the vectors in the tangent space at $(0, 0, \bar{z})$ for $\bar{z} = 66.302 \dots$
converging to the unstable eigenvectors $(1, 0, 0)^T$ and $(0, 1, 0)^T$ at $O$, respectively.
We define two initial points, $u_1$ and $u_2$,
created by perturbing the transient orbit in the directions of $v_1$ and $v_2$;
$u_1=(0, 0, \bar{z})^T + s v_1, u_2=(0, 0, \bar{z})^T + s v_2$
where $s = 10^{-12}$ is a small amplitude of the perturbation.
Fig. \ref{perturbation-of-first-vector} and Fig. \ref{perturbation-of-second-vector}
shows the orbits starting from $u_1$ and $u_2$, respectively.
The two orbits in the figure go toward xy plane near z-axis in an initial interval,
and the direction changes drastically in a neighborhood of $O$.
The projections of the orbits to xy plane show that the orbits go to the expected directions;
the directions parallel to the corresponding eigenvectors of the linearized system at $O$.

\section{Conclusion}

We have proposed a calculation algorithm
to analyze tangent spaces of transient orbits converging to an equilibrium,
which is a direct application of the method by Ginelli et al.
to calculate covariant Lyapunov vectors numerically.
The application gives vectors in tangent spaces of a transient orbit that
converge to eigenvectors of the linearized system at the equilibrium.

To test our algorithm, we analyzed a simple model of three-dimensional ordinary differential equation.
We also showed that our suggested vectors for a transient orbit converging to an equilibrium are useful
to give initial points whose orbit passes
near a direction of a desired eigenvector having an unstable eigenvalue.
Actually, we succeeded in demonstrating that giving two initial points near the stable manifold so that
their orbit moving in the directions of eigenvectors in a neighborhood of the equilibrium.

Transient orbits play important roles in many dynamical problems.
In particular, pattern formation problems such as self-replicating patterns and
collision problems of localized structures may have some relationship with our methods.
The authors have been analyzing such problems using the present method,
and their detailed results will be given elsewhere.

\bibliographystyle{plain}
\bibliography{references_en}

\end{document}